\documentclass{ijclclp}

\title{The Turán number of the balanced double star $S_{n-1,n-1}$ in the hypercube $Q_n$}
\author{Dan-dan Liu,  
        Shou-jun Xu*}
\affiliation{* Corresponding author.

\noindent{School of Mathematics and Statistics, Gansu Center for Applied Mathematics, Lanzhou University, Lanzhou, Gansu 730000, China}}
\email{\texttt{liudd2023@lzu.edu.cn, shjxu@lzu.edu.cn}}
\usepackage{amsthm}

\usepackage{tipx}  
\DeclareUnicodeCharacter{030B}{\doubleacute}  
\usepackage{doi}
\usepackage[numbers]{natbib}

\begin{document}
\maketitle

\thispagestyle{firstpage}
\begin{abstract}
The $n$-dimensional hypercube $Q_n$ is a graph with vertex set $\{0,1\}^n$ such that there is an edge between two vertices if and only if they differ in exactly one coordinate. Let $H$ be a graph, and a graph is called $H$-free if it does not contain $H$ as a subgraph. Given a graph $H$, the Turán number of $H$ in $Q_n$, denoted by $ex(Q_n, H)$, is the maximum number of edges of a subgraph of
$Q_n$ that is $H$-free. A double star $S_{k,l}$ is the graph obtained by taking an edge $uv$ and joining $u$ with $k$ vertices,
$v$ with $l$ vertices which are different from the $k$ vertices. We say a double star is a balanced double star if $k = l$. Currently, the Turán number of the balanced star $S_{n,n}$ is has been studied in the planar graphs. In the hypercubes, the researchers look for the maximum number of edges of the graphs that are $C_k$-free. However, the Turán number of the double star in the hypercube remains unexplored. Building upon prior research, we initiate the first study on the Turán number of the balanced double star in the hypercube. In this paper, we give the exact value of the Turán number of the balanced double star $S_{n-1,n-1}$ in the hypercube $Q_n$, which is $ 2^{n-3}\times(4n- 3)$ for all $n \geq 3$.

\end{abstract}

\section{Introduction }
 
 Let $G=(V(G), E(G))$ be a graph, where $V(G)$ and $E(G)$ are the vertex set and edge set of the graph $G$, respectively. We use $N_G(v)$ to denote the set of vertices of $G$ adjacent to $v$, abbreviated without ambiguity as $N(v)$. The degree of a vertex $v$ is the number of vertices of $N(v)$, denoted by $d_G(v)$, abbreviated without ambiguity as $d(x)$ i.e. $|N(v)|$. Let $N_G[v]=N(v) \cup \{v\} $, abbreviated without ambiguity as $N[v]$. Let $v(G)$, $e(G)$, $\delta(G)$ and $\Delta(G)$ denote the number of
vertices, the number of edges, the minimum degree and the maximum degree of $G$, respectively. We use $E[S,T]$ to denote the edge set between $S$ and $T$, and let $e[S,T]$ = $|E[S,T]|$, where $S,T$ are subsets of $V(G)$.

A graph is called $H$-free if it does not contain the graph $H$ as a subgraph. For graphs $G$ and $H$, let $ex(G, H)$ denote the maximum number of edges in a subgraph
of $G$ which does not contain a copy of $H$. The $n$-dimensional hypercube $Q_n$ is a graph with vertex set $\{0,1\}^n$ such that there is an edge between two vertices if and only if they differ in exactly one coordinate. A double star $S_{k,l}$ is the graph obtained by taking an edge $uv$ and joining $u$ with $k$ vertices,
$v$ with $l$ vertices which are different from the $k$ vertices. We say a double star is a balanced double star if $k = l$.

 Significant progress has been made in the study of Turán-type problems, where the base graph $G$ is the $n$-dimensional hypercube $Q_n$. Erdős 
\cite{P11111} initiated the study of $ex(Q_n, H)$, the maximum number of edges in an $H$-free subgraph of $Q_n$ in the special case when $H$ is an even cycle. He asked how many edges a $C_4$-free subgraph of the hypercube can contain and conjectured that the answer is
($\tfrac{1}{2}+ o(1))e(Q_n)$ and offered $100$ for a solution. The best known upper bound, due to Baber \cite{Baber}, is approximately
$0.60318 · e(Q_n)$. Erdős \cite{P11111} also raised the extremal question for even cycles. Chung \cite{Chung} obtained $\tfrac{ex(Q_n,C_{4k})} {e(Q_n)})$ $\rightarrow$ 0 for every $k$ $\leq$ 2 , i.e. cycles whose lengths are divisible by 4, starting from lengths 8 are harder to avoid than the four-cycles. Later Conder \cite{Conder} improved the lower bound to $\tfrac{1}{3}$$e(Q_n)$ by defining a 3-coloring of the
edges of the $n$-cube such that every color class is $C_6$-free. And Alexandr \cite{Alexandr} $et$ $al$.
proved that for any $n$ $\in$ $N$, $ex(Q_n, C_{10})\textgreater  0.024 · e(Q_n)$.
 
 The double star is an important structure in graph theory. However, the special structure of the double star has not been considered in the planar graphs or hypercubes before. In 2022, Győri, Martin, Paulos and Xiao \cite{2021arXiv211010515G} studied the topic for double stars as the forbidden graph. A double star $S_{k,l}$ is the graph obtained by taking an edge $uv$ and joining $u$ with $k$ vertices,
$v$ with $l$ vertices which are different from the $k$ vertices. We say a double star is a balanced double star if $k = l$. A $k$-$l$ edge refers to an edge whose endpoints have degrees $k$ and $l$, respectively. Győri, Martin, Paulos and Xiao established the exact value for $S_{2,2}$ and $S_{2,3}$, and derived the upper bounds of $S_{2,4}$, $S_{2,5}$, $S_{3,3}$, $S_{3,4}$. In 2024, Xu \cite{XU2024326,2024arXiv240901016X, xu2024planarturannumberbalanced} $et$ $al$. improved the upper bounds of $S_{2,4}$, $S_{2,5}$ and $S_{3,3}$. In 2025, Xu \cite{XU} $et$ $al$. improved the upper bounds of $S_{3,4}$,  Liu \cite{liu} gave the upper bounds of $S_{3,5}$.

In this paper, we first investigate the Turán number of the balanced double star $S_{n-1,n-1}$ in the hypercube $Q_n$, and determine its exact value.

\textbf{Theorem 1.1.} Let $Q_n$ be the $n$-dimensional hypercube on $2^n$ vertices, Then

$$
ex(Q_n, S_{n-1,n-1}) = \left\{
\begin{array}{ll}
0,    &if\quad n=1 \\
2,    & if\quad n=2 \\
 2^{n-3}\times(4n- 3), &if\quad n \geq 3 \\
\end{array} 
\right.
$$

The remainder of this paper is organized as follows. Essential preliminaries are provided in Section 2. Section 3 establishes the proof of Theorem 1.1, with detailed analysis of its implications.

\section{Preliminaries}
\textbf{Claim 2.1. }
For any $n$-dimensional hypercube $Q_n$ with $n \geq 3$, let $G_n$  be a subgraph of $Q_n$ that is $S_{n-1,n-1}$-free and has the maximum number of edges. If $\delta(G_n)$ $\leq$ $n-2$, then $G_n$ can be transformed into a subgraph $G_n'$ with the same number of edges such that $\delta(G_n')$ = $n-1$ and $G_n'$ remains $S_{n-1,n-1}$-free.

\noindent\textbf{Proof.} Here, $G_n$ being $S_{n-1,n-1}$-free implies that no two vertices of degree $n$ are adjacent in $G_n$. 
Assume that $\delta(G_n)$= $t$ $\leq$ $n-2$. Then there exists a vertex $v$ $\in$ $G_n$ with 
$d_{G_n}(v)$= $t$ $\leq$ $n-2$. Since $G_n$ is a subgraph of 
$Q_n$, $v$ is missing at least two edges 
$vp_1$ and $vq_1$ from $Q_n$.

Suppose that we attempt to add the edge $vp_1$ to $G_n$. If $d_{G_n}(p_1)$$\leq$ $n-2$, adding $vp_1$ would result in $d_{G_n}(v)$ $\leq$ $n-1$ and 
$d_{G_n}(p_1)$$\leq$ $n-1$, which appears permissible. However, this leads to a contradiction: If there exists an edge could be added to 
$G_n$ without creating an $S_{n-1,n-1}$, then $G_n$ would not be an $S_{n-1,n-1}$-free subgraph of $Q_n$ with the maximum number of edges.

Now $d_{G_n}(p_1)$ = $n-1$. Consider that the neighborhood 
$N_{G_n}(p_1)$. If every vertex $s$ $\in$ $N_{G_n}(p_1)$ satisfies $d_{G_n}(s)$ $\leq$ $n-1$, then we can add the edge $vp_1$ without any constraints, since $d_{G_n}(p_1)$ = $n$ after adding the edge $vp_1$. Hence, there exists a vertex $s_1$ of $N_{G_n}(p_1)$ with $d_{G_n}(s)$ = $n$. Then we can delete the edge $p_1s_1$ and add the edge $vp_1$ to construct the $G_n''$. i.e. $G_n''$=($G_n$ $\backslash \{p_1s_1\}$) $\cup \{vp_1\}$. Then $\delta(G_n'')$ = $t+1$. 

Through $n-1-t$ such operations, we can transform $G_n$ into $G_n'$. i.e. $G_n'$ = ($G_n$ $\backslash$ ($\{p_1s_1\}$ $\cup$ $\{p_2s_2\}$ $\cup$ $\cdots$ $\cup$ $\{p_{n-1-t}s_{n-1-t}\}$)) $\cup$ ($\{vp_1\}$ $\cup$ 
$\{vp_2\}$ $\cup$ $\cdots$ $\cup$ $\{vp_{n-1-t}\}$).$\hfill\square$

Suppose that $G_n$ is a $S_{n-1,n-1}$-free subgraph of $Q_n$ with the maximum number of edges . If $\delta(G_n)$ $\leq$ $n-2$, we can apply the aforementioned method to transform $G_n$ into another $S_{n-1,n-1}$-free subgraph $G_n'$ $\subseteq$ $Q_n$ with $\delta(G_n')$ = $n-1$, while maintaining the same number of edges ($e(G_n)=e(G_n')$). Consequently, we may restrict our analysis to maximum-edge
$S_{n-1,n-1}$-free subgraph $G_n$ $\subseteq$ $Q_n$ satisfying $\delta(G_n)$ = $n-1$.
\vspace{3mm}

\noindent\textbf{Claim 2.2. } For any $n$-dimensional hypercube $Q_n$ with $n \geq 3$, there exist at least two distinct $S_{n-1,n-1}$-free subgraphs $G_n$ and $G_n'$ of $Q_n$ after labelling, each containing $2^{n-3}\times (4n-3)$ edges. Furthermore, under any vertex labelling scheme, no vertex has the same degree $n$ in both $G_n$ and $G_n'$.

\noindent\textbf{Proof.} (1) Base case. For $n=3$, there exist two distinct $S_{2,2}$-free subgraphs $G_3$ and $G_3'$ of $Q_3$, each attaining the maximal edge count $2^{3-3}\times(4\times3-3)$ = 9, such that no vertex of degree 3 is common to both  $G_3$ and $G_3'$. Their structural differences are illustrated in the Figure 1.

\begin{figure}[h]
    \centering
    \includegraphics[width=0.7\textwidth]{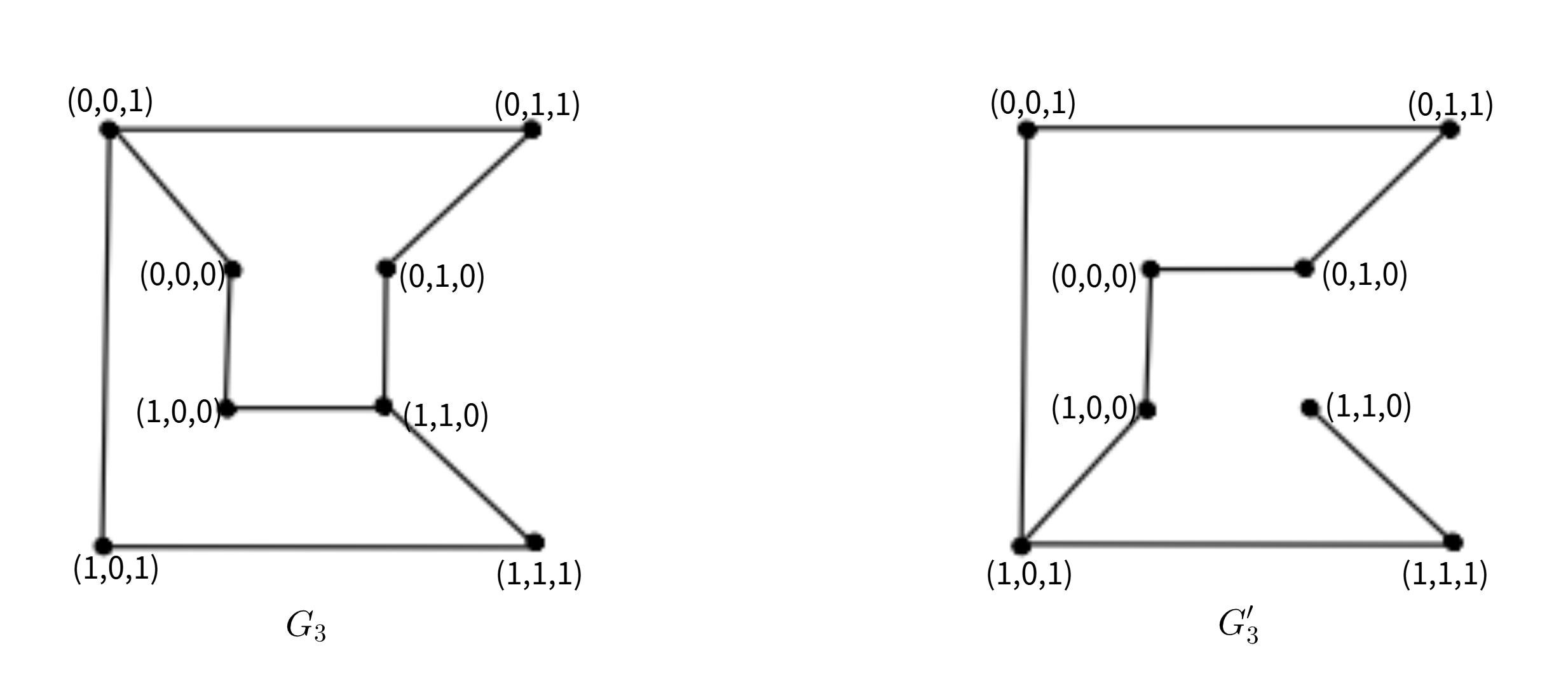}
    \caption{$G_3$ and $G_3'$.}
    \label{} 
\end{figure}

(2) Inductive Step. Suppose that the assertion holds for all hypercubes $Q_t$ with $t \leq n$; that is there exist two distinct $S_{t-1,t-1}$-free subgraphs $G_t$ and $G_t'$ of $Q_t$ after labelling, each attaining the maximal edge count $2^{t-3}\times(4t-3)$ , such that no vertex of degree $t$ is common to both  $G_t$ and $G_t'$.

We now consider $Q_{n+1}$. By the Cartesian product structure $Q_{n+1}$ = $Q_n$$\times$$K_2$, we extend $G_n$ and $G_n'$ $\subseteq$ $Q_n$ as follows (see Figure 2):

(i) Constructing $G_{n+1}$: Assign the 
$(n+1)$-th coordinate of all vertices in $G_n$  to 0 and those in $G_n'$ to 1. Their Cartesian product forms $G_{n+1}$ $\subseteq$ $Q_{n+1}$.

(ii) Constructing $G_{n+1}'$: Assign the 
$(n+1)$-th coordinate of all vertices in $G_n$  to 1 and those in $G_n'$ to 0. Their Cartesian product forms $G_{n+1}'$ $\subseteq$ $Q_{n+1}$.

\begin{figure}[h]
    \centering
    \includegraphics[width=0.9\textwidth]{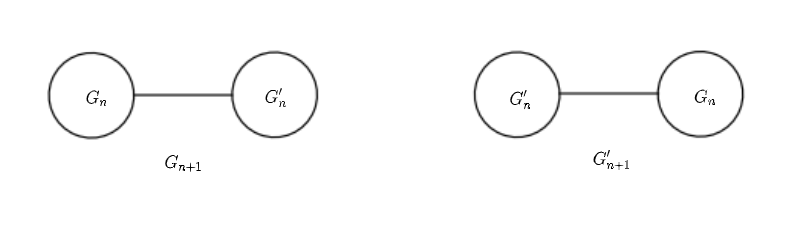}
    \caption{$G_{n+1}$ and $G_{n+1}'$.}
    \label{Fig:57354} 
\end{figure}

\noindent $e(G_n)$ = $e(G_n')$ = $2^{n-3}\times (4n-3)$. Moreover, after labelling, no vertex of degree $n$ appears in both $G_n$ and $G_n'$. And $e[G_n,G_n']$ = $2^n$. We can add there edges without any constraints. So $e(G_{n+1})$ = $e(G_{n+1}')$ = 
($2^{n-3}\times (4n-3)$) $\times 2$+$2^n$ = $2^{n-2}\times (4(n+1)-3)$. Then there exist two distinct subgraphs $G_{n+1}$ and $G_{n+1}'$ of $Q_{n+1}$ after labelling that satisfy the conditions.$\hfill\square$

\vspace{3mm}      
 
\noindent\textbf{Claim 2.3. } For any $n$-dimensional  hypercube $Q_n$ with $n \geq 3$, there does not exist an $S_{n-1,n-1}$-free subgraph $G_n$ $\subseteq$ $Q_n$ with $e(G_n)$ $\geq$ $2^{n-3}\times (4n-3)+1$. i.e. any $S_{n-1,n-1}$-free subgraph of $Q_n$ must remove at least $3\times 2^{n-3}$ edges from $Q_n$, then $e(G_n)$ $\leq$ 
 $2^{n-3}\times (4n-3)$.

\noindent\textbf{Proof.} 
(1) Base case. For $n=3$, Suppose $G_3$ $\subseteq$ $Q_3$ is $S_{2,2}$-free with $e(G_3)$ $\geq$ $10$. By the degree-sum formula, $G_3$ must contain at least four vertices of degree 3. Let $S$= $\{u_1,u_2,u_3,u_4\}$ denote these vertices of degree $3$. Let $S_i$ denote the set of the vertices of $Q_n$ with exactly $i$ coordinates equal to 1, so that $|S_i|$= $\binom{n}{i}$ and each vertex in $S_i$ is adjacent to  $i$ vertices in $S_{i-1}$ and $n-i$ vertices in $S_{i+1}$. Partition the vertices of 
$Q_3$ into four sets $S_0$, $S_1$, $S_2$, $S_3$. By symmetry, we have $S_0$ $\cong$ $S_3$ and $S_1$ $\cong$ $S_2$.

We assert that $S_0\cup S_3$ contains at most one vertex from $S$. For contradiction, assume one vertex $u_1$ $\in$ $S\cap S_0$ and another vertex $u_2$ $\in$ $S\cap S_3$. The remaining two vertices $\{u_3,u_4\}$ $\in$ $S$ must belong to $S_1\cup  S_2$. By the degree constraint $d_{G_3}(u_3)$ = $d_{G_3}(u_4)$ = 3, both $u_3$ and $u_4$ must each have edges connecting to 
$u_1$ or $u_2$. This implies: if either $u_3$ or $u_4$ belongs to $S_1$, then it must be adjacent to  $u_1$ $\in$ $S_0$; if either $u_3$ or $u_4$ belongs to $S_2$, then it must be adjacent to  $u_2$ $\in$ $S_3$.
Then at least two vertices of degree $3$ are adjacent, forming a balanced double star $S_{2,2}$, a contradiction. 

 Next, we discuss the location of vertices in $S$ (see Figure 3).

\begin{figure}[h]
    \centering
    \includegraphics[width=1.0\textwidth]{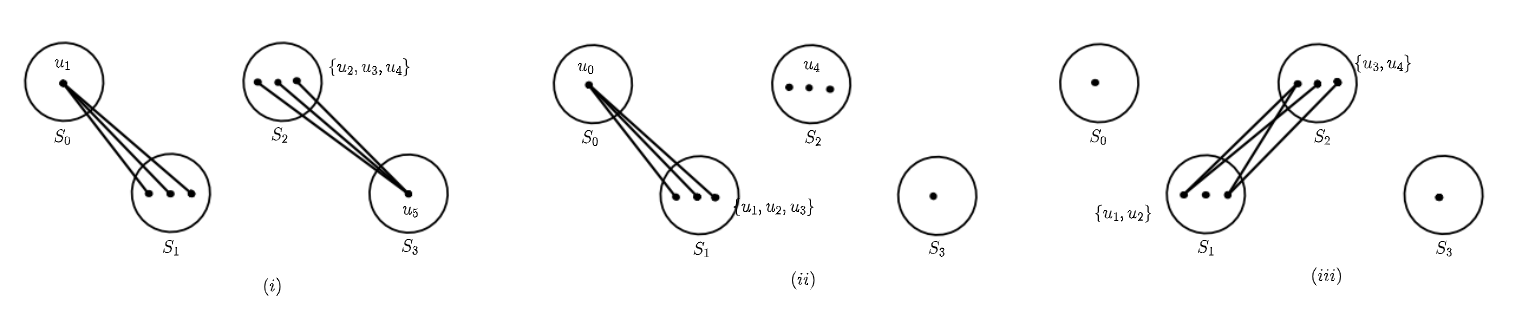}
    \caption{ the location of vertices in $S$.}
    \label{Fig:57354} 
\end{figure}

Case 1: There is exactly one vertex in $S \cap(S_0\cup S_3)$. Assume that $u_1$ $\in$ $S\cap S_0$, as illustrated in the Figure 3(i).

No vertex in $\{u_2,u_3,u_4\}$ can belong to $S_1$. Otherwise, such a vertex would be adjacent to $u_1$, forming a balanced double star $S_{2,2}$, which is contradictory. Furthermore, $\{u_2,u_3,u_4\}$ $\nsubseteq$ $S_3$. Therefore, $\{u_2,u_3,u_4\}$ must lie in $S_2$. i.e. $\{u_2,u_3,u_4\}$ = $S_2$. Let $u_5$ denote the vertex in $S_3$. Since all neighbors of $u_5$ in $G_3$ have degree 3, it follows that $d_{G_3}(u_5)$ = 3. This configuration necessarily contains a balanced double star $S_{2,2}$ in $G_3$, a contradiction. 

Case 2: There are no vertices in $S \cap(S_0\cup S_3)$. There are two possibilities.

Case 2.1: Exactly one vertex in $S\cap S_2$, and the remaining three vertices are in $S_1$, as illustrated in the Figure 3(ii).

Assume that $\{u_1,u_2,u_3\}$ $\subseteq$ $S_1$ and $u_4$ $\subseteq$ $S_2$. Let $u_0$ denote the vertex in $S_0$. 
 Since all neighbors of $u_0$ in $G_3$ have degree 3, we have $d_{G_3}(u_0)$ = 3. This configuration necessarily contains a balanced double star $S_{2,2}$ in $G_3$, a contradiction. 

Case 2.2: Two vertices in $S\cap S_2$ and two in $S\cap S_1$, as illustrated in the Figure 3(ii). We hereby establish the validity of the conclusion under these conditions through mathematical induction.

 Here, $u_1$ and $u_2$ share exactly one common neighbor in $S_2$, forcing at least one adjacency between $\{u_1,u_2\}$ and $\{u_3,u_4\}$. This creates two vertices of degree 3 are adjacent, forming a balanced double star $S_{2,2}$ in $G_3$, a contradiction.

In conclusion, for $n=3$, if $G_3$ is an $S_{2,2}$-free subgraph of $Q_3$, then $e(G_3)$ $\leq$ $9$. 

(2) Inductive Step. Assume that  for $t \leq n$ every $S_{t-1,t-1}$-free subgraph $G_t$ of $Q_t$ satisfies $e(G_t)$ $\leq$ $2^{t-3}\times (4t-3)$.

Consider that the hypercube  $Q_{n+1}$ = $Q_n$$\times$$K_2$. Suppose that there exists an $S_{n,n}$-free subgraph $G_{n+1}$ $\subseteq$ $Q_{n+1}$ with $e(G_{n+1})$ $\geq$ $2^{n-2}\times(4(n+1)-3)+1$.
\vspace{3mm}

\noindent\textbf{Subclaim 1. }The balanced double star $S_{n,n}$ exist in any edge-deleted subgraph $G_{n+1}$ $\subseteq$ $Q_{n+1}$ with $e(G_{n+1})$ $\geq$ $2^{n-2}\times(4(n+1)-3)+1$, regardless of how edges are removed from $e[Q_n,Q_n]$.

\noindent\textbf{Proof.} 
(1)\textbf {No edges deletions in $e[Q_n,Q_n]$.}

At least one copy of $Q_n$  (left or right) removes 
at most 3 × $2^{n-3}-1$ edges. By induction hypothesis, this copy contains an $S_{n-1,n-1}$. With intact intercube edges, this extends to an $S_{n,n}$ in $G_{n+1}$$\subseteq$ $Q_{n+1}$, contradicting the $S_{n,n}$-free assumption.

(2)\textbf{Edge deletions in $e[Q_n,Q_n]$.}

Suppose that deleting at most $r$ intercube edges in $e[Q_n,Q_n]$, there is no $S_{n,n}$-free subgraph $G_{n+1}$$\subseteq$ $Q_{n+1}$ with $e(G_{n+1})$ $\geq$ $2^{n-2}\times(4(n+1)-3)+1$. 

\begin{figure}[h]
    \centering
    \includegraphics[width=1.0\textwidth]{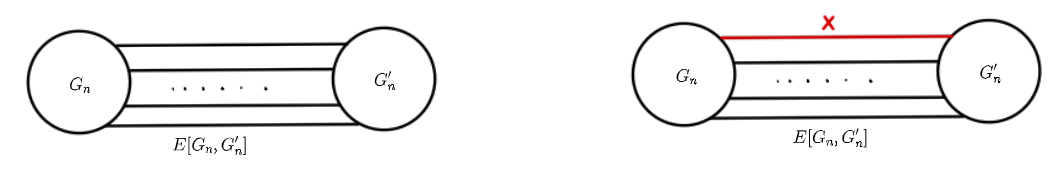}
    \caption{ the edge deletions in $e[Q_n,Q_n]$.}
    \label{Fig:57354} 
\end{figure}

If we delete one additional edge from $e[Q_n,Q_n]$, and reintroduce it into the left subgraph $G_n$ $\subseteq$ $Q_n$ (see Figure 4), the minimum degree condition $\delta(G_n)$ $ \geq$ $n-1$ partitions the vertices of $G_n$ into three classes (as shown in Figure 5):

(i) $R_1$: Vertices with full degree $n$ in $G_n$, i.e. each vertex $s$ of $G_n$ is of degree $n$ (each intercube edges to $G_n'$ are removed).

(ii) $R_2$: Vertices of degree $n$ in $G_{n+1}$ retaining all intercube connections to the right subgraph $G_n'$.

(iii) $R_3$: Vertices having degree $n+1$ in $G_{n+1}$ while preserving intercube edges to $G_n'$.

\begin{figure}[h]
    \centering
    \includegraphics[width=0.7\textwidth]{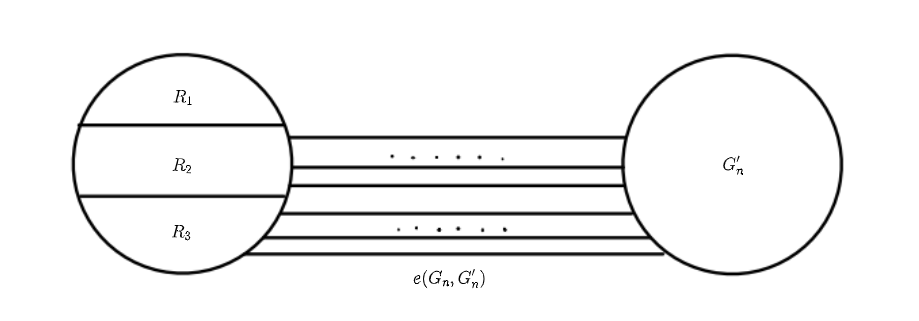}
    \caption{$G_{n+1}$ $\subseteq$ $Q_{n+1}$ with $e(G_{n+1})$ $\geq$ $2^{n-2}\times(4(n+1)-3)+1$ and $e[Q_n,Q_n]$=$2^n-(r+1)$. }
    \label{Fig:57354} 
\end{figure}

Thus, the additional edge can only be added to $R_2$. However, no matter how this edge is added, two vertices of degree $n+1$ in $R_2$ must become adjacent, thereby forming the balanced double star $S_{n,n}$, which is a contradiction.

This contradiction refutes the assumptions, thereby validating the original conclusion.

Consequently, Claim 2.3 holds as asserted.$\hfill\square$

\vspace{3mm}

\section{Proof of Theorem 1.1.}
\noindent\textbf{Theorem 1.1.} Let $Q_n$ be a $n$-dimensional hypercube on $2^n$ vertices, Then

$$
ex(Q_n, S_{n-1,n-1}) = \left\{
\begin{array}{ll}
0,    & n=1 \\
2,    & n=2 \\
n \times 2^{n-1} - 3 \times 2^{n-3}, & n \geq 3 \\
\end{array} 
\right.
$$

\noindent\textbf{Proof.}
For $n=1, 2$, the statement holds trivially.
Let $G_n$ be an $S_{n-1,n-1}$-free subgraph of $Q_n$. 

By Claim 2.1, let $G_n$ be a $S_{n-1,n-1}$-free subgraph of $Q_n$ with the maximum number of edges. If $\delta(G_n)$ $\leq$ $n-2$, we can transform $G_n$ into another $S_{n-1,n-1}$-free subgraph $G_n'$ $\subseteq$ $Q_n$ with $\delta(G_n')$ = $n-1$, while maintaining the same number of edges ($e(G_n)=e(G_n')$). So we just need to consider that $\delta(G_n)$ = $n-1$.

By Claim 2.2, for any $n$-dimensional hypercube $Q_n$ with $n \geq 3$, there exist at least two distinct $S_{n-1,n-1}$-free subgraphs $G_n$ and $G_n'$ of $Q_n$ after labelling, each containing $2^{n-3}\times (4n-3)$ edges.

Conversely, by Claim 2.3, we konw that any $S_{n-1,n-1}$-free subgraph of $Q_n$ must remove at least $3\times 2^{n-3}$ edges from $Q_n$, then $e(G_n)$ $\leq$ 
 $2^{n-3}\times (4n-3)$.

Combining these results yields $e(G_n)$=$2^{n-3}\times (4n-3)$, thereby establishing Theorem 1.1.
$\hfill\square$

\section*{Funding}
This work was funded in part by Natural Science Foundation of China (Grant No. 12071194).

\section*{Data availability}
No data was used for the research described in the article.

\bibliographystyle{plain}
\bibliography{main}
\end{document}